\documentclass[11pt]{article}
\usepackage{amsfonts,amscd,amssymb}
\newtheorem{definition}{Definition}[section]

\newtheorem{proposition}[definition]{Proposition}
\newtheorem{corollary}[definition]{Corollary}
\newtheorem{remark}[definition]{Remark}

\newtheorem{theorem}[definition]{Theorem}
\newtheorem{example}[definition]{Example}

\def\rawo\lonra{\longrightarrow}

\def\ot{\otimes}

\newcommand{\eqref}[1]{(\ref{eq:#1})}

\newenvironment{proof}{{\it Proof.}}{\hfill $ \square $ \vskip 4mm}
\begin{document}
\title{Generalized (anti) Yetter-Drinfeld modules as components of 
a braided T-category}
\author
{Florin Panaite\thanks{Research partially supported by the programme 
CERES of the Romanian Ministry of Education and Research, contract 
no. 4-147/2004.}\\
Institute of Mathematics of the 
Romanian Academy\\ 
PO-Box 1-764, RO-014700 Bucharest, Romania\\
e-mail: Florin.Panaite@imar.ro\\
\and 
Mihai D. Staic\thanks{Permanent address: Institute of Mathematics of the  
Romanian Academy, 
PO-Box 1-764, RO-014700 Bucharest, Romania.}\\
SUNY at Buffalo\\
Amherst, NY 14260-2900, USA\\
e-mail: mdstaic@buffalo.edu}
\date{}
\maketitle
%%%%%%%%%%%%%%%%%%%%%%%%%%%%%%%%%%%%%%%%%%%%%%%%%%%%%%%%%
\begin{abstract}
If $H$ is a Hopf algebra with bijective antipode and $\alpha , \beta $$\;\in  
Aut_{Hopf}(H)$, we introduce a category $_H{\cal YD}^H(\alpha , \beta )$,  
generalizing both Yetter-Drinfeld modules and anti-Yetter-Drinfeld modules. 
We construct a braided T-category ${\cal YD}(H)$  
having all the categories $_H{\cal YD}^H(\alpha , \beta )$  
as components, which if $H$ is finite dimensional coincides with the 
representations of a certain quasitriangular 
T-coalgebra $DT(H)$ that we construct. We also   
prove that if $(\alpha , \beta )$ admits a so-called  pair in   
involution, then $_H{\cal YD}^H(\alpha , \beta )$ is isomorphic to the 
category of usual Yetter-Drinfeld modules $_H{\cal YD}^H$.   
\end{abstract}
%%%%%%%%%%%%%%%%%%%%%%%%%%%%%%%%%%%%%%%%%%%%%%%%%%%
\section*{Introduction}
%%%%%%%%%%%%%%%%%%%%%%%%%%%%%%%%%%%%%%%%%%%%%%%%%%%%
Let $H$ be a Hopf algebra with bijective antipode $S$ and 
$\alpha , \beta $$\; 
\in Aut_{Hopf}(H)$. We introduce the concept of an  
$(\alpha , \beta )$-{\it Yetter-Drinfeld  
module}, as being a left $H$-module right $H$-comodule $M$ with the  
following compatibility condition: 
\begin{eqnarray*}
&&(h\cdot m)_{(0)}\ot (h\cdot m)_{(1)}=h_2\cdot m_{(0)}\ot 
\beta (h_3)m_{(1)}\alpha (S^{-1}(h_1)). 
\end{eqnarray*}
This concept is a generalization of three kinds of objects appeared in the 
literature. Namely, for $\alpha =\beta =id_H$, we obtain the usual 
Yetter-Drinfeld modules; for $\alpha =S^2$, $\beta =id_H$, we obtain the 
so-called anti-Yetter-Drinfeld modules, introduced in \cite{hajac1}, 
\cite{hajac2}, \cite{jara} as coefficients for the cyclic cohomology of 
Hopf algebras defined by Connes and Moscovici in \cite{connes1},  
\cite{connes2}; finally, an $(id_H, \beta )$-Yetter-Drinfeld module is a 
generalization of the object $H_{\beta }$ defined in \cite{cvoz}, 
which has the property that, if $H$ is finite dimensional, then the map 
$\beta \mapsto End(H_{\beta })$ gives a group anti-homomorphism from 
$Aut_{Hopf}(H)$ to the Brauer group of $H$. \\
It is natural to expect that $(\alpha , \beta )$-Yetter-Drinfeld modules 
have some properties resembling the ones of the three kinds of objects 
we mentioned. We will see some of these properties in this paper (others 
will be given in a subsequent one), namely the ones directed to our main 
aim here, which is the following: if we denote by 
$_H{\cal YD}^H(\alpha , \beta )$ the category of $(\alpha , \beta )$
-Yetter-Drinfeld modules and we define ${\cal YD}(H)$ as the disjoint 
union of all these categories, then we can organize ${\cal YD}(H)$ as a 
braided T-category (or braided crossed group-category, in the original  
terminology of Turaev, see \cite{turaev}) over the group 
$G=Aut_{Hopf}(H)\times Aut_{Hopf}(H)$ with multiplication 
$(\alpha , \beta )*(\gamma , \delta )=(\alpha \gamma , \delta \gamma ^{-1}
\beta \gamma )$. We also prove that the subcategory 
${\cal YD}(H)_{fd}$ consisting of finite dimensional objects has left 
and right dualities, and that, if $H$ is finite dimensional, then 
${\cal YD}(H)$ coincides with the representations of a certain quasitriangular 
T-coalgebra $DT(H)$ that we construct.\\
Our second aim is to prove that, if $\alpha , \beta $$\;\in   
Aut_{Hopf}(H)$ such that there exists a so-called {\it pair in involution} 
$(f, g)$ corresponding to $(\alpha , \beta )$, then 
$_H{\cal YD}^H(\alpha , \beta )$ is isomorphic to $_H{\cal YD}^H$. This 
result is independent on the theory concerning ${\cal YD}(H)$, but we 
can give it a very short proof using the results obtained during the 
construction of ${\cal YD}(H)$.  
%%%%%%%%%%%%%%%%%%%%%%%%%%%%%%%%%%%%%%%%%%%%%%%%%%%%%%%%%%%%%%
\section{Preliminaries}\label{sec1}
%%%%%%%%%%%%%%%%%%%%%%%%%%%%%%%%%%%%%%%%%%%%%%%%%%%%%%%%%%%%%%%
We work over a ground field $k$. All algebras, linear spaces,  
etc. will be over $k$; unadorned $\ot $ means $\ot_k$. Unless otherwise 
stated, $H$ will denote a Hopf algebra with bijective antipode $S$. We will 
use the versions of Sweedler's sigma notation: $\Delta (h)=h_1\ot h_2$ or 
$\Delta (h)=h_{(1)}\ot h_{(2)}$.  
For unexplained concepts and notation about Hopf algebras we refer to 
\cite{k}, \cite{m}, \cite{mon}, \cite{sw}.  
By $\alpha , \beta, \gamma ...$ we will usually denote Hopf automorphisms 
of $H$.\\[2mm]
Let $A$ be 
an $H$-bicomodule algebra, with comodule structures $A\rightarrow
A\ot H$, $a\mapsto a_{<0>}\ot a_{<1>}$ and $A\rightarrow H\ot A$,
$a\mapsto a_{[-1]}\ot a_{[0]}$, and denote, for $a\in A$,
\begin{eqnarray*}
&&a_{\{-1\}}\ot a_{\{0\}}\ot a_{\{1\}}=a_{<0>_{[-1]}}\ot
a_{<0>_{[0]}}\ot a_{<1>}=a_{[-1]}\ot a_{[0]_{<0>}}\ot
a_{[0]_{<1>}}
\end{eqnarray*}
as an element in $H\ot A\ot H$. We can consider the 
{\it Yetter-Drinfeld datum} $(H, A, H)$ as in
\cite{cmz} (the second $H$ is regarded as an $H$-bimodule
coalgebra), and the Yetter-Drinfeld category $_A {\cal YD}(H)^H$,
whose objects are $k$-modules $M$ endowed with a left $A$-action
(denoted by $a\ot m\mapsto a\cdot m$) and a right $H$-coaction
(denoted by $m\mapsto m_{(0)}\ot m_{(1)}$) satisfying the equivalent 
compatibility conditions
\begin{eqnarray}
&&(a\cdot m)_{(0)}\ot (a\cdot m)_{(1)}=a_{\{0\}}\cdot m_{(0)}
\ot a_{\{1\}}m_{(1)}S^{-1}(a_{\{-1\}}),\\ \label{ydd1}
&&a_{<0>}\cdot m_{(0)}\ot a_{<1>}m_{(1)}=(a_{[0]}\cdot
m)_{(0)}\ot (a_{[0]}\cdot m)_{(1)}a_{[-1]}, \label{ydd2}
\end{eqnarray}
for all $a\in A$ and $m\in M$.\\
Recall now from 
\cite{hn1} the construction of the (left) {\it diagonal crossed product} 
$H^*\bowtie A$, which is an associative algebra constructed on 
$H^*\ot A$, with 
multiplication given by
\begin{eqnarray}
&&(p\bowtie a)(q\bowtie b)=p(a_{\{-1\}}\rightharpoonup q\leftharpoonup
S^{-1}(a_{\{1\}}))\bowtie a_{\{0\}}b,
\end{eqnarray}
for all $a, b\in A$ and $p, q\in H^*$, 
and with unit $\varepsilon _H\bowtie 1_A$. Here 
$\rightharpoonup $ and $\leftharpoonup $ are
the regular actions of $H$ on $H^*$ given by 
$(h\rightharpoonup p)(l)=p(lh)$ and $(p\leftharpoonup
h)(l)=p(hl)$ for all $h, l\in H$ and $p\in H^*$.\\
If $H$ is finite dimensional, we can consider 
the Drinfeld double $D(H)$, which is a quasitriangular Hopf algebra realized 
on $H^*\ot H$; its coalgebra structure is $H^{* cop}\ot H$ and the algebra 
structure is just $H^*\bowtie H$, that is 
\begin{eqnarray}
&&(p\bowtie h)(q\bowtie l)=p(h_1\rightharpoonup q\leftharpoonup 
S^{-1}(h_3))\bowtie h_2l,
\end{eqnarray}
for all $p, q\in H^*$ and $h, l\in H$.\\
The diagonal crossed product  
$H^*\bowtie A$ becomes a $D(H)$-bicomodule algebra, with
structures
\begin{eqnarray*}
&&H^*\bowtie A\rightarrow (H^*\bowtie A)\ot D(H), \ p\bowtie a\mapsto
(p_2\bowtie a_{<0>})\ot (p_1\ot a_{<1>}), \\
&&H^*\bowtie A\rightarrow D(H)\ot (H^*\bowtie A), \ p\bowtie
a\mapsto (p_2\ot a_{[-1]})\ot (p_1\bowtie a_{[0]}),
\end{eqnarray*}
for all $p\in H^*$ and $a\in A$, see \cite{hn1}.\\
In the case when $H$ is finite dimensional, by results in 
\cite{bpvo}, \cite{cmz} it follows that the category $_A{\cal YD}(H)^H$ is 
isomorphic to the category $_{H^*\bowtie A}{\cal M}$ of left modules over 
$H^*\bowtie A$.   
%%%%%%%%%%%%%%%%%%%%%%%%%%%%%%%%%%%%%%%%%%%%%%%%%
\section{$(\alpha , \beta )$-Yetter-Drinfeld modules}\label{sec2}
%%%%%%%%%%%%%%%%%%%%%%%%%%%%%%%%%%%%%%%%%%%%%%%%
\setcounter{equation}{0}
%%%%%%%%%%%%%%%%%%%%%%%%%%%%%%%%%%%%%%%%%%%%%%%%%%%%%%%%
\begin{definition} Let $\alpha , \beta $$\;\in Aut_{Hopf}(H)$.   
An $(\alpha , \beta )$-Yetter-Drinfeld module over $H$ 
is a vector space $M$, such that $M$ is a left $H$-module (with notation 
$h\ot m\mapsto h\cdot m$) and a right $H$-comodule (with notation 
$M\rightarrow M\ot H$, $m\mapsto m_{(0)}\ot m_{(1)}$) with the following  
compatibility condition:
\begin{eqnarray}
&&(h\cdot m)_{(0)}\ot (h\cdot m)_{(1)}=h_2\cdot m_{(0)}\ot 
\beta (h_3)m_{(1)}\alpha (S^{-1}(h_1)), \label{ab}
\end{eqnarray}
for all $h\in H$ and $m\in M$. We denote by $_H{\cal YD}^H(\alpha , \beta )$  
the category of $(\alpha , \beta )$-Yetter-Drinfeld modules, morphisms 
being the $H$-linear $H$-colinear maps.
\end{definition}
\begin{remark}
As for usual Yetter-Drinfeld modules, one can see that (\ref{ab}) is 
equivalent to 
\begin{eqnarray}
&&h_1\cdot m_{(0)}\ot \beta (h_2)m_{(1)}=(h_2\cdot m)_{(0)}\ot 
(h_2\cdot m)_{(1)}\alpha (h_1). \label{abn}
\end{eqnarray}
\end{remark}
\begin{example}
For $\alpha =\beta =id_H$, we have $_H{\cal YD}^H(id, id)=$$_H{\cal YD}^H$,  
the usual category of (left-right) Yetter-Drinfeld modules. 
\end{example} 
\begin{example}
For $\alpha =S^2$, $\beta =id_H$, the compatibility condition (\ref{ab}) 
becomes 
\begin{eqnarray}
&&(h\cdot m)_{(0)}\ot (h\cdot m)_{(1)}=h_2\cdot m_{(0)}\ot 
h_3m_{(1)}S(h_1), \label{ayd}
\end{eqnarray}
hence $_H{\cal YD}^H(S^2, id)$ is the category of anti-Yetter-Drinfeld 
modules defined in \cite{hajac1}, \cite{hajac2}, \cite{jara}.
\end{example} 
\begin{example} \label{cvz}
For $\beta $$\;\in Aut_{Hopf}(H)$, define $H_{\beta }$ as in \cite{cvoz}, 
that is $H_{\beta }=H$, with regular right $H$-comodule structure and 
left $H$-module structure given by $h\cdot h'=\beta (h_2)h'S^{-1}(h_1)$, for 
all $h, h'\in H$. It was noticed in \cite{cvoz} that $H_{\beta }$ satisfies a  
certain compatibility condition, which actually says that 
$H_{\beta }\in $$\;_H{\cal YD}^H(id, \beta )$. More generally, if   
$\alpha , \beta $$\;\in Aut_{Hopf}(H)$, define $H_{\alpha , \beta }$ as 
follows: $H_{\alpha , \beta }=H$, with regular right $H$-comodule structure 
and left $H$-module structure given by 
$h\cdot h'=\beta (h_2)h'\alpha (S^{-1}(h_1))$, for $h, h'\in H$. Then one 
can check that $H_{\alpha , \beta }\in $$\;_H{\cal YD}^H(\alpha , \beta )$. 
\end{example}    
\begin{example}
Take $l$ an integer and define $\alpha _l=S^{2l}$$\;\in Aut_{Hopf}(H)$. The 
compatibility in $_H{\cal YD}^H(S^{2l}, id)$ becomes 
\begin{eqnarray}
&&(h\cdot m)_{(0)}\ot (h\cdot m)_{(1)}=h_2\cdot m_{(0)}\ot 
h_3m_{(1)}S^{2l-1}(h_1). 
\end{eqnarray}
An object in $_H{\cal YD}^H(S^{2l}, id)$ will be called an 
$l-{\cal YD}$-module. Hence, a $0-{\cal YD}$-module is a Yetter-Drinfeld 
module and a $1-{\cal YD}$-module is an anti-Yetter-Drinfeld module. The 
right-left version of $l-{\cal YD}$-modules has been introduced in 
\cite{doru}. 
\end{example}
\begin{example}\label{mp1}
Let $\alpha $, $\beta$$\;\in Aut_{Hopf}(H)$ and assume that there exist 
an algebra map $f:H\rightarrow k$ and a group-like element $g\in H$ such that
\begin{eqnarray} 
&&\alpha (h)=g^{-1}f(h_1)\beta(h_2)f(S(h_3))g, \;\;\;\forall \;\;h\in H. 
\label{pi2} 
\end{eqnarray}
Then one can   
check that $k\in $$\;_H{\cal YD}^H(\alpha , \beta)$, with structures 
$\;h\cdot 1=f(h)$ and $1\mapsto 1\ot g$. More generally, if $V$ is any 
vector space, then $V\in $$\;_H{\cal YD}^H(\alpha , \beta)$, with 
structures $h\cdot v=f(h)v$ and $v\mapsto v\ot g$, for all $h\in H$ and 
$v\in V$. 
\end{example}
\begin{definition}
If $\alpha $, $\beta$$\;\in Aut_{Hopf}(H)$ such that there exist $f, g$ as in 
Example \ref{mp1}, we will say that $(f, g)$ is a pair in involution 
corresponding to $(\alpha,\beta)$ (in analogy with the concept of modular pair 
in involution due to Connes and Moscovici)  
and the $(\alpha , \beta)$-Yetter-Drinfeld modules 
$k$ and $V$ constructed in Example \ref{mp1} will be denoted by $_fk^g$ and 
respectively $_fV^g$. 
\end{definition}
As an example, if $\alpha $$\;\in Aut_{Hopf}(H)$, then  
$(\varepsilon, 1)$ is a pair in involution corresponding to 
$(\alpha , \alpha )$.\\[2mm]   
Let $\alpha , \beta $$\;\in Aut_{Hopf}(H)$. We define an $H$-bicomodule  
algebra $H(\alpha , \beta )$ as follows: $H(\alpha , \beta )=H$ as algebra, 
with comodule structures
\begin{eqnarray*}
&&H(\alpha , \beta )\rightarrow H\ot H(\alpha , \beta ), \;\;\;\;
h\mapsto h_{[-1]}\ot h_{[0]}=\alpha (h_1)\ot h_2, \\
&&H(\alpha , \beta )\rightarrow H(\alpha , \beta )\ot H, \;\;\;\;
h\mapsto h_{<0>}\ot h_{<1>}=h_1\ot \beta (h_2).
\end{eqnarray*}
Then we can consider the Yetter-Drinfeld datum $(H, H(\alpha , \beta ), H)$ 
and the Yetter-Drinfeld modules over it, 
$_{H(\alpha , \beta )}{\cal YD}(H)^H$.
\begin{proposition}
$_H{\cal YD}^H(\alpha , \beta )=$$\;_{H(\alpha , \beta )}{\cal YD}(H)^H$.
\end{proposition}
\begin{proof}
It is easy to see that the compatibility conditions for the two categories 
are the same.  
\end{proof}   
In particular, the category of anti-Yetter-Drinfeld modules coincides 
with $_{H(S^2, id)}{\cal YD}(H)^H$, which improves the remark in 
\cite{hajac1} that anti-Yetter-Drinfeld modules are entwined modules. \\[2mm]
Consider now the diagonal crossed product 
$A(\alpha , \beta )=H^*\bowtie H(\alpha , \beta )$, whose multiplication is 
\begin{eqnarray}
&&(p\bowtie h)(q\bowtie l)=p(\alpha (h_1)\rightharpoonup q\leftharpoonup  
S^{-1}(\beta (h_3)))\bowtie h_2l, 
\end{eqnarray}
for all $p, q\in H^*$ and $h, l\in H$. For $\alpha =\beta =id$ we get 
$A(id, id)=D(H)$; for $\alpha =S^2$ and $\beta =id$, the multiplication 
in $A(S^2, id)$ is 
\begin{eqnarray}
&&(p\bowtie h)(q\bowtie l)=p(S^2(h_1)\rightharpoonup q\leftharpoonup  
S^{-1}(h_3))\bowtie h_2l, 
\end{eqnarray}
hence $A(S^2, id)$ coincides with the algebra $A(H)$ defined in 
\cite{hajac1}.\\[2mm]
Assume now that $H$ is finite dimensional; then $A(\alpha , \beta )$ 
becomes a $D(H)$-bicomodule algebra, with structures 
\begin{eqnarray*}
&&H^*\bowtie H(\alpha , \beta )
\rightarrow (H^*\bowtie H(\alpha , \beta ))\ot D(H), \;\;\; 
p\bowtie h\mapsto 
(p_2\bowtie h_1)\ot (p_1\ot \beta (h_2)), \\
&&H^*\bowtie H(\alpha , \beta )
\rightarrow D(H)\ot (H^*\bowtie H(\alpha , \beta )), \;\;\; 
p\bowtie 
h\mapsto (p_2\ot \alpha (h_1))\ot (p_1\bowtie h_2). 
\end{eqnarray*}
In particular, $A(H)$ becomes a $D(H)$-bicomodule algebra, improving the 
remark in \cite{hajac1} that $A(H)$ is a right $D(H)$-comodule algebra. 
Since $H$ is finite dimensional, we have an isomorphism of categories 
$_{H(\alpha , \beta )}{\cal YD}(H)^H\simeq  
$$\;_{H^*\bowtie H(\alpha , \beta )}{\cal M}$, hence    
$_H{\cal YD}^H(\alpha , \beta )\simeq  
$$\;_{H^*\bowtie H(\alpha , \beta )}{\cal M}$ (for $\alpha =S^2$, 
$\beta =id$    
we recover the result in \cite{hajac1} that the category of 
anti-Yetter-Drinfeld modules is isomorphic to $_{A(H)}{\cal M}$). The 
correspondence is given as follows. If 
$M\in $$\;_H{\cal YD}^H(\alpha , \beta )$, then 
$M\in $$\;_{H^*\bowtie H(\alpha , \beta )}{\cal M}$ with structure 
\begin{eqnarray*}
&&(p\bowtie h)\cdot m=p((h\cdot m)_{(1)})(h\cdot m)_{(0)}.
\end{eqnarray*}
Conversely, if $M\in $$\;_{H^*\bowtie H(\alpha , \beta )}{\cal M}$, 
then $M\in $$\;_H{\cal YD}^H(\alpha , \beta )$ with structures 
\begin{eqnarray*}
&&h\cdot m=(\varepsilon \bowtie h)\cdot m, \\
&&m\mapsto m_{(0)}\ot m_{(1)}=(e^i\bowtie 1)\cdot m\ot e_i,
\end{eqnarray*}
where $\{e_i\}$, $\{e^i\}$ are dual bases in $H$ and $H^*$.
%%%%%%%%%%%%%%%%%%%%%%%%%%%%%%%%%%%%%%%%%%%%%%%%%%%%%%%%%%%%%%%%%%%
\section{A braided T-category ${\cal YD}(H)$}\label{sec3}
%%%%%%%%%%%%%%%%%%%%%%%%%%%%%%%%%%%%%%%%%%%%%%%%%%%%%%%%%%%%%%%%%
\setcounter{equation}{0}
%%%%%%%%%%%%%%%%%%%%%%%%%%%%%%%%%%%%%%%%%%%%%%%%%%%%%%%%%%%%%%
Let $\alpha , \beta $$\;\in Aut_{Hopf}(H)$ and consider the objects 
$H_{\alpha }, H_{\beta }$ as in Example \ref{cvz}. In \cite{cvoz} was 
considered the object $M=H_{\alpha }\ot H_{\beta }$, with the 
following structures:
\begin{eqnarray*}
&&h\cdot (x\ot y)=h_1\cdot x\ot \alpha (h_2)\cdot y,\\
&&x\ot y\mapsto (x_1\ot y_1)\ot y_2x_2,
\end{eqnarray*}
for all $h, x, y\in H$, where by $\cdot $ we denoted both the actions of 
$H$ on $H_{\alpha }$ and $H_{\beta }$ given as in Example \ref{cvz}.  
Then it was noticed in \cite{cvoz} that $M$ satisfies a compatibility 
condition which says that 
$M\in $$\;_H{\cal YD}^H(id, \beta \alpha )$.\\
On the other hand, it was noticed in \cite{hajac1} that the tensor product 
between an anti-Yetter-Drinfeld module and a Yetter-Drinfeld module becomes 
an anti-Yetter-Drinfeld module. \\  
The next result can be seen as a generalization of both these facts.  
\begin{proposition} \label{unu}
If $M\in $$\;_H{\cal YD}^H(\alpha , \beta )$,  
$N\in $$\;_H{\cal YD}^H(\gamma , \delta )$, then  
$M\ot N\in $$\;_H{\cal YD}^H(\alpha \gamma , 
\delta \gamma ^{-1}\beta \gamma )$, with structures:
\begin{eqnarray*}
&&h\cdot (m\ot n)=\gamma (h_1)\cdot m\ot 
\gamma ^{-1}\beta \gamma (h_2)\cdot n, \\
&&m\ot n\mapsto (m\ot n)_{(0)}\ot (m\ot n)_{(1)}=(m_{(0)}\ot n_{(0)})\ot 
n_{(1)}m_{(1)}.
\end{eqnarray*}
\end{proposition}
\begin{proof}
Obviously $M\ot N$ is a left $H$-module and a right $H$-comodule. We check 
now the compatibility condition. We compute: \\[2mm]
${\;\;\;\;}$$(h\cdot (m\ot n))_{(0)}\ot (h\cdot (m\ot n))_{(1)}$
\begin{eqnarray*}
&&=(\gamma (h_1)\cdot m \ot \gamma ^{-1}\beta \gamma (h_2)\cdot n)_{(0)}\ot  
(\gamma (h_1)\cdot m\ot \gamma ^{-1}\beta \gamma (h_2)\cdot n)_{(1)}\\
&&=((\gamma (h_1)\cdot m)_{(0)}\ot (\gamma ^{-1}\beta \gamma (h_2)
\cdot n)_{(0)})\ot (\gamma ^{-1}\beta \gamma (h_2)\cdot n)_{(1)}
(\gamma (h_1)\cdot m)_{(1)}\\
&&=(\gamma (h_1)_2\cdot m_{(0)}\ot \gamma ^{-1}\beta \gamma (h_2)_2\cdot 
n_{(0)})\ot \\ 
&&\;\;\;\;\;\;\;\;\;\;\;\;\;\;\ot \; 
\delta (\gamma ^{-1}\beta \gamma (h_2)_3)n_{(1)}\gamma   
(S^{-1}(\gamma ^{-1}\beta \gamma (h_2)_1))\beta (\gamma (h_1)_3)m_{(1)}
\alpha (S^{-1}(\gamma (h_1)_1))\\
&&=(\gamma (h_2)\cdot m_{(0)}\ot \gamma ^{-1}\beta \gamma (h_5)\cdot 
n_{(0)})\ot \delta \gamma ^{-1}\beta \gamma (h_6)n_{(1)}
\beta \gamma (S^{-1}(h_4))\beta \gamma (h_3)m_{(1)}
\alpha (S^{-1}(\gamma (h_1)))\\
&&=(\gamma (h_2)\cdot m_{(0)}\ot \gamma ^{-1}\beta \gamma (h_3)\cdot 
n_{(0)})\ot \delta \gamma ^{-1}\beta \gamma (h_4)n_{(1)}m_{(1)}
\alpha \gamma (S^{-1}(h_1))\\
&&=h_2\cdot (m_{(0)}\ot n_{(0)})\ot \delta \gamma ^{-1}\beta \gamma (h_3)
n_{(1)}m_{(1)}\alpha \gamma (S^{-1}(h_1))\\
&&=h_2\cdot (m\ot n)_{(0)}\ot \delta \gamma ^{-1}\beta \gamma (h_3)
(m\ot n)_{(1)}\alpha \gamma (S^{-1}(h_1)),
\end{eqnarray*}
that is $M\ot N\in $$\;_H{\cal YD}^H(\alpha \gamma ,  
\delta \gamma ^{-1}\beta \gamma )$. 
\end{proof}
Note that, if $M\in $$\;_H{\cal YD}^H(\alpha , \beta )$,    
$N\in $$\;_H{\cal YD}^H(\gamma , \delta )$ and 
$P\in $$\;_H{\cal YD}^H(\mu , \nu )$, then $(M\ot N)\ot P=
M\ot (N\ot P)$ as objects in 
$_H{\cal YD}^H(\alpha \gamma \mu , \nu \mu ^{-1}\delta \gamma ^{-1}
\beta \gamma \mu )$.\\[2mm]
Denote $G=Aut_{Hopf}(H)\times Aut_{Hopf}(H)$, a group with multiplication 
\begin{eqnarray}
&&(\alpha , \beta )*(\gamma , \delta )=(\alpha \gamma , \delta \gamma ^{-1}
\beta \gamma )
\end{eqnarray}
(the unit is $(id, id)$ and $(\alpha , \beta )^{-1}=(\alpha ^{-1}, 
\alpha \beta ^{-1}\alpha ^{-1}))$.
\begin{proposition}
Let $N\in $$\;_H{\cal YD}^H(\gamma , \delta )$ and $(\alpha , \beta )\in G$. 
Define $^{(\alpha , \beta )}N=N$ as vector space, with structures 
\begin{eqnarray*}
&&h\rightharpoonup n=\gamma ^{-1}\beta \gamma \alpha ^{-1}(h)\cdot n, \\
&&n\mapsto n_{<0>}\ot n_{<1>}=n_{(0)}\ot \alpha \beta ^{-1}(n_{(1)}). 
\end{eqnarray*}
Then $^{(\alpha , \beta )}N\in $$\;_H{\cal YD}^H
(\alpha \gamma \alpha ^{-1}, \alpha \beta ^{-1}\delta \gamma ^{-1}\beta 
\gamma \alpha ^{-1})=$$\;_H{\cal YD}^H((\alpha , \beta )*(\gamma , \delta )* 
(\alpha , \beta )^{-1})$. 
\end{proposition}
\begin{proof}
Obviously $^{(\alpha , \beta )}N$ is a left $H$-module and right 
$H$-comodule, so we check the compatibility condition. We compute:\\[2mm]
${\;\;\;\;\;}$$(h\rightharpoonup n)_{<0>}\ot (h\rightharpoonup n)_{<1>}$
\begin{eqnarray*}
&&=(\gamma ^{-1}\beta \gamma \alpha ^{-1}(h)\cdot n)_{(0)}\ot 
\alpha \beta ^{-1}((\gamma ^{-1}\beta \gamma \alpha ^{-1}(h)\cdot n)_{(1)})\\
&&=\gamma ^{-1}\beta \gamma \alpha ^{-1}(h_2)\cdot n_{(0)}\ot 
\alpha \beta ^{-1}(\delta \gamma ^{-1}\beta \gamma \alpha ^{-1}(h_3)
n_{(1)}\gamma \gamma ^{-1}\beta \gamma \alpha ^{-1}(S^{-1}(h_1)))\\
&&=\gamma ^{-1}\beta \gamma \alpha ^{-1}(h_2)\cdot n_{(0)}\ot 
\alpha \beta ^{-1}\delta \gamma ^{-1}\beta \gamma \alpha ^{-1}(h_3)
\alpha \beta ^{-1}(n_{(1)})\alpha \gamma \alpha ^{-1}(S^{-1}(h_1))\\
&&=h_2\rightharpoonup n_{(0)}\ot \alpha \beta ^{-1}\delta \gamma ^{-1}\beta 
\gamma \alpha ^{-1}(h_3)n_{<1>}\alpha \gamma \alpha ^{-1}(S^{-1}(h_1)),
\end{eqnarray*}
that is $^{(\alpha , \beta )}N\in $$\;_H{\cal YD}^H 
(\alpha \gamma \alpha ^{-1}, \alpha \beta ^{-1}\delta \gamma ^{-1}\beta 
\gamma \alpha ^{-1})$.   
\end{proof} 
\begin{remark} \label{gigi}
Let $M\in $$\;_H{\cal YD}^H(\alpha , \beta )$,   
$N\in $$\;_H{\cal YD}^H(\gamma , \delta )$ and $(\mu , \nu )\in G$. Then we  
have  
\begin{eqnarray*}
&&^{(\alpha , \beta )*(\mu , \nu )}N=\;^{(\alpha , \beta )}(^{(\mu , \nu )}N) 
\end{eqnarray*}
as objects in $_H{\cal YD}^H 
(\alpha \mu \gamma \mu ^{-1}\alpha ^{-1}, 
\alpha \beta ^{-1}\mu \nu ^{-1}\delta \gamma ^{-1}\nu \mu ^{-1}
\beta \mu \gamma \mu ^{-1}\alpha ^{-1})$, and 
\begin{eqnarray*}
&&^{(\mu , \nu )}(M\ot N)=\;^{(\mu , \nu )}M\ot \;^{(\mu , \nu )}N 
\end{eqnarray*}
as objects in $_H{\cal YD}^H(\mu \alpha \gamma \mu ^{-1}, 
\mu \nu ^{-1}\delta \gamma ^{-1}\beta \alpha ^{-1}\nu \alpha \gamma 
\mu ^{-1})$.
\end{remark}
\begin{proposition}
Let $M\in $$\;_H{\cal YD}^H(\alpha , \beta )$ and    
$N\in $$\;_H{\cal YD}^H(\gamma , \delta )$. Define $^MN=\;
^{(\alpha , \beta )}N$ as object in $_H{\cal YD}^H  
((\alpha , \beta )*(\gamma , \delta )*(\alpha , \beta )^{-1})$. Define 
the map
\begin{eqnarray*}
&&c_{M, N}:M\ot N\rightarrow \;^MN\ot M,\;\;\;c_{M, N}(m\ot n)=
n_{(0)}\ot \beta ^{-1}(n_{(1)})\cdot m.
\end{eqnarray*}
Then $c_{M, N}$ is $H$-linear $H$-colinear and satisfies the conditions 
(for $P\in $$\;_H{\cal YD}^H(\mu , \nu )$):
\begin{eqnarray}
&&c_{M\ot N, P}=(c_{M, \;^NP}\ot id_N)\circ (id_M\ot c_{N, P}), \label{br1}\\
&&c_{M, N\ot P}=(id_{\;^MN}\ot c_{M, P})\circ (c_{M, N}\ot id_P).\label{br2}
\end{eqnarray}
Moreover, if $M\in $$\;_H{\cal YD}^H(\alpha , \beta )$,     
$N\in $$\;_H{\cal YD}^H(\gamma , \delta )$ and $(\mu , \nu )\in G$, then 
$c_{\;^{(\mu ,\nu )}M, \;^{(\mu , \nu )}N}=c_{M, N}$. 
\end{proposition}
\begin{proof}
We prove that $c_{M, N}$ is $H$-linear. We compute:\\[2mm]
${\;\;\;\;\;}$
$c_{M, N}(h\cdot (m\ot n))$
\begin{eqnarray*}
&&=c_{M, N}(\gamma (h_1)\cdot m\ot  
\gamma ^{-1}\beta \gamma (h_2)\cdot n)\\
&&=(\gamma ^{-1}\beta \gamma (h_2)\cdot n)_{(0)}\ot \beta ^{-1}
((\gamma ^{-1}\beta \gamma (h_2)\cdot n)_{(1)})\gamma (h_1)\cdot m \\
&&=\gamma ^{-1}\beta \gamma (h_2)_2\cdot n_{(0)}\ot \beta ^{-1}(\delta  
(\gamma ^{-1}\beta \gamma (h_2)_3)n_{(1)}\gamma (S^{-1}((\gamma ^{-1}\beta 
\gamma (h_2))_1)))\gamma (h_1)\cdot m \\
&&=\gamma ^{-1}\beta \gamma (h_3)\cdot n_{(0)}\ot \beta ^{-1}\delta  
\gamma ^{-1}\beta \gamma (h_4)\beta ^{-1}(n_{(1)})\gamma (S^{-1}(h_2))
\gamma (h_1)\cdot m\\
&&=\gamma ^{-1}\beta \gamma (h_1)\cdot n_{(0)}\ot \beta ^{-1}\delta  
\gamma ^{-1}\beta \gamma (h_2)\beta ^{-1}(n_{(1)})\cdot m,
\end{eqnarray*}
\begin{eqnarray*}  
h\cdot c_{M, N}(m\ot n)&=&h\cdot (n_{(0)}\ot \beta ^{-1}(n_{(1)})\cdot m)\\
&=&\alpha (h_1)\rightharpoonup n_{(0)}\ot \alpha ^{-1}\alpha \beta ^{-1}
\delta \gamma ^{-1}\beta \gamma \alpha ^{-1}\alpha (h_2)\beta ^{-1}(n_{(1)})
\cdot m\\
&=&\gamma ^{-1}\beta \gamma \alpha ^{-1}\alpha (h_1)\cdot n_{(0)}\ot  
\beta ^{-1}\delta \gamma ^{-1}\beta \gamma (h_2)\beta ^{-1}(n_{(1)})\cdot m\\
&=&\gamma ^{-1}\beta \gamma (h_1)\cdot n_{(0)}\ot \beta ^{-1}\delta  
\gamma ^{-1}\beta \gamma (h_2)\beta ^{-1}(n_{(1)})\cdot m,
\end{eqnarray*}
so the two terms are equal. The fact that $c_{M, N}$ is $H$-colinear is 
similar and left to the reader. We prove now (\ref{br1}). First note that, 
due to Remark \ref{gigi}, we have $^M(^NP)=$$\;^{M\ot N}P$ and $^M(N\ot P)=$$ 
\;^MN\ot \;^MP$. We compute: 
\begin{eqnarray*}
(c_{M, \;^NP}\ot id_N)\circ (id_M\ot c_{N, P})(m\ot n\ot p)&=&
c_{M, \;^NP}(m\ot p_{(0)})\ot \delta ^{-1}(p_{(1)})\cdot n\\
&=&p_{(0)_{<0>}}\ot \beta ^{-1}(p_{(0)_{<1>}})\cdot m\ot \delta ^{-1}
(p_{(1)})\cdot n\\
&=&p_{(0)_{(0)}}\ot \beta ^{-1}\gamma \delta ^{-1}
(p_{(0)_{(1)}})\cdot m\ot \delta ^{-1}(p_{(1)})\cdot n\\
&=&p_{(0)}\ot \beta ^{-1}\gamma \delta ^{-1}
(p_{(1)_1})\cdot m\ot \delta ^{-1}(p_{(1)_2})\cdot n, 
\end{eqnarray*}   
\begin{eqnarray*}
c_{M\ot N, P}(m\ot n\ot p)&=&p_{(0)}\ot \gamma ^{-1}\beta ^{-1}\gamma 
\delta ^{-1}(p_{(1)})\cdot (m\ot n)\\
&=&p_{(0)}\ot \gamma \gamma ^{-1}\beta ^{-1}\gamma  
\delta ^{-1}(p_{(1)_1})\cdot m\ot \gamma ^{-1}\beta \gamma \gamma ^{-1}
\beta ^{-1}\gamma \delta ^{-1}(p_{(1)_2})\cdot n\\
&=&p_{(0)}\ot \beta ^{-1}\gamma \delta ^{-1}
(p_{(1)_1})\cdot m\ot \delta ^{-1}(p_{(1)_2})\cdot n,
\end{eqnarray*}
and we are done. The proof of (\ref{br2}) is easier and left to the reader, 
and similarly the last statement of the Proposition.
\end{proof} 
Note that $c_{M, N}$ is bijective with inverse $c_{M, N}^{-1}(n\ot m)=
\beta ^{-1}(S(n_{(1)}))\cdot m\ot n_{(0)}$. \\[2mm]
We are ready now to introduce the desired braided T-category (we use 
terminology as in \cite{zunino}; for the subject of Turaev categories, see 
also the original paper of Turaev \cite{turaev} and \cite{stef}, 
\cite{virelizier}).\\
Define ${\cal YD}(H)$ as the disjoint union of all $_H{\cal YD}^H(\alpha , 
\beta )$, with $(\alpha , \beta )\in G$ (hence the component of the unit 
is just $_H{\cal YD}^H$). If we endow ${\cal YD}(H)$ with tensor product 
as in Proposition \ref{unu}, then it becomes a strict monoidal category 
with unit $k$ as object in $_H{\cal YD}^H$ (with trivial structures). \\
The group homomorphism $\varphi :G\rightarrow aut ({\cal YD}(H))$, 
$(\alpha , \beta )\mapsto \varphi _{(\alpha , \beta )}$, is given on 
components as 
\begin{eqnarray*}
&&\varphi _{(\alpha , \beta )}:\;_H{\cal YD}^H(\gamma , \delta )
\rightarrow \;_H{\cal YD}^H((\alpha , \beta )*(\gamma , \delta )
*(\alpha , \beta )^{-1}), \;\;\;\varphi _{(\alpha , \beta )}(N)=\;
^{(\alpha , \beta )}N, 
\end{eqnarray*}
and the functor $\varphi _{(\alpha , \beta )}$ acts as identity on morphisms. 
The braiding in ${\cal YD}(H)$ is given by the family 
$\{c_{M, N}\}$. 
As a consequence of the above results, we obtain:
\begin{theorem}
${\cal YD}(H)$ is a braided T-category over $G$. 
\end{theorem}
We consider now the problem of existence of left and right dualities. 
\begin{proposition}\label{abp} 
Let $M\in $$\;_H{\cal YD}^H(\alpha , \beta )$ and assume 
that $M$ is finite dimensional. Then $M^*=Hom(M,k)$ becomes an object in   
$_H{\cal YD}^H(\alpha^{-1} , \alpha\beta^{-1}\alpha^{-1})$,   
with $(h\cdot f) (m)=f((\beta^{-1}\alpha^{-1}S(h))\cdot m)$ and    
$f_{(0)}(m)\ot f_{(1)}=f(m_{(0)})\ot S^{-1}(m_{(1)})$. 
Moreover, the maps $b_M:k\to M\ot M^*$, $b_M(1)=\sum_ie_i\ot e^i$ (where 
$\{e_i\}$ and $\{e^i\}$ are dual bases in $M$ and $M^*$)  
and $d_M:M^*\ot M\to k$, $d_M(f\ot m)=f(m)$, are morphisms in 
$_H{\cal YD}^H$ and we have  
$(id_M\ot d_M)(b_M\ot id_M)=id_M$ and $(d_M\ot id_{M^*})(id_{M^*}\ot b_M)=
id_{M^*}$.
\end{proposition}
\begin{proof}
We first prove that      
$M^*$ is indeed an object in $\;_H{\cal YD}^H(\alpha^{-1} , \alpha\beta^{-1} 
\alpha^{-1} )$. We compute:
\begin{eqnarray*}
(h\cdot f)_{(0)}(m)\ot (h\cdot f)_{(1)}&=
&(h\cdot f)(m_{(0)})\ot S^{-1}(m_{(1)})\\
&=&f((\beta^{-1}\alpha^{-1}S)(h)\cdot m_{(0)})\ot S^{-1}(m_{(1)}), 
\end{eqnarray*}
${\;\;\;\;}$$(h_{(2)}\cdot f_{(0)})(m)\ot 
(\alpha \beta^{-1} \alpha^{-1})(h_{(3)}) f_{(1)} 
(\alpha^{-1}S^{-1})(h_{(1)})$
\begin{eqnarray*}
%&&=(h_{(2)}\cdot f_{(0)})(m)\ot (\alpha\beta^{-1}\alpha^{-1}(h_{(3)}))f_{(1)}
%(\alpha^{-1}S^{-1})(h_{(1))})\\
&&=f_{(0)}(\beta^{-1}\alpha^{-1}S(h_{(2)})\cdot m)
\ot (\alpha\beta^{-1}\alpha^{-1})
(h_{(3)})f_{(1)}(\alpha^{-1}S^{-1}(h_{(1)}))\\
&&=f(((\beta^{-1}\alpha^{-1}S)(h_{(2)})\cdot m)_{(0)})\ot  
(\alpha\beta^{-1}\alpha^{-1}(h_{(3)}))\\
&&\;\;\;\;\;\;\;\;\;\;\;\;\;\;\;
S^{-1}(((\beta^{-1}\alpha^{-1}S)(h_{(2)})\cdot m)_{(1)})
(\alpha^{-1}S^{-1}(h_{(1)}))\\
&&=f((\beta^{-1}\alpha^{-1}S)(h_{(3)})\cdot m_{(0)})\ot 
(\alpha\beta^{-1}\alpha^{-1}(h_{(5)}))\\
&&\;\;\;\;\;\;\;\;\;\;\;\;\;\;\;S^{-1}((\alpha^{-1}S)(h_{(2)})m_{(1)}
(\alpha\beta^{-1}\alpha^{-1})(h_{(4)}))(\alpha^{-1}S^{-1}(h_{(1)}))\\
&&=f((\beta^{-1}\alpha^{-1}S)(h_{(3)})\cdot m_{(0)})\ot 
(\alpha\beta^{-1}\alpha^{-1})
(h_{(5)}S^{-1}(h_{(4)}))S^{-1}(m_{(1)})\alpha^{-1}(h_{(2)}S^{-1}(h_{(1)}))\\
&&=f((\beta^{-1}\alpha^{-1}S)(h)\cdot m_{(0)})\ot S^{-1}(m_{(1)}), 
\end{eqnarray*}
which means that
\begin{eqnarray*}
&&(h\cdot f)_{(0)}\ot (h\cdot f)_{(1)}=(h_{(2)}\cdot f_{(0)})\ot  
(\alpha \beta^{-1} \alpha^{-1})(h_{(3)}) f_{(1)} (\alpha^{-1}S^{-1})(h_{(1)}), 
\;\;q.e.d.
\end{eqnarray*}
On $k$ we have the trivial module and comodule structure, and  
with these $k\in$$\;_H{\cal YD}^H$.  
We want to prove that $b_M$ and $d_M$ are $H$-module maps. We compute: 
\begin{eqnarray*}
(h\cdot b_M(1))(m)&=&(h\cdot(\sum_i e_i\ot e^i))(m)\\
&=&\sum_i \alpha^{-1}(h_{(1)})\cdot e_i\ot 
((\alpha\beta\alpha^{-1})(h_{(2)})\cdot e^i)(m)\\
&=&\alpha^{-1}(h_{(1)})\cdot e_i\ot e^i((\beta^{-1}\alpha^{-1}
S\alpha\beta\alpha^{-1})(h_{(2)})\cdot m)\\
&=&\sum_i \alpha^{-1}(h_{(1)})\cdot e_i\ot e^i((\alpha^{-1}S)(h_{(2)})
\cdot m)\\
&=&\alpha^{-1}(h_{(1)}S(h_{(2)}))\cdot m\\
&=&\varepsilon(h)\sum_i e_i\ot e^i(m)\\
&=&(\varepsilon(h)b_M(1))(m),
\end{eqnarray*}
\begin{eqnarray*}
d_M(h\cdot (f\ot m))&=&d_M(\alpha(h_{(1)})\cdot f\ot \beta^{-1}
(h_{(2)})\cdot m)\\
&=&(\alpha(h_{(1)})\cdot f)(\beta^{-1}(h_{(2)})\cdot m)\\
&=&f((\beta^{-1}\alpha^{-1}S\alpha(h_{(1)}))\beta^{-1}(h_{(2)})\cdot m)\\
&=&f(\beta^{-1}(S(h_{(1)})h_{(2)})\cdot m)\\
&=&\varepsilon(h)d_M(f\ot m).
\end{eqnarray*}
Also they are $H$-comodule maps:
\begin{eqnarray*}
((b_M(1))_{(0)}\ot (b_M(1))_{(1)})(m)&=
&\sum_i (e_i)_{(0)} \ot (e^i)_{(0)}(m) 
\ot (e^i)_{(1)}(e_i)_{(1)}\\
&=&\sum_i (e_i)_{(0)}\ot (e^i)(m_{(0)})\ot S^{-1}(m_{(1)})(e_i)_{(1)}\\
&=&m_{(0)}\ot S^{-1}(m_{(1)_2})m_{(1)_1}\\
&=&(b_M(1)\ot 1)(m),
\end{eqnarray*}
\begin{eqnarray*}
d_M((f\ot m)_{(0)})\ot (f\ot m)_{(1)}&=& f_{(0)}(m_{(0)})\ot m_{(1)}f_{(1)}\\
&=&f(m_{(0)})\ot m_{(1)_2}S^{-1}(m_{(1)_1})\\
&=&d_M(f\ot m)\ot 1.
\end{eqnarray*}
Finally, the last two identities $(id_M\ot d_M)(b_M\ot id_M)=id_M$  
and $(d_M\ot id_{M^*})(id_{M^*}\ot b_M)=id_{M^*}$ are trivial.
\end{proof}
Similarly, one can prove:
\begin{proposition}\label{rightdual} 
Let $M\in $$\;_H{\cal YD}^H(\alpha , \beta )$ and  
assume that $M$ is finite dimensional.  
Then $^*M=Hom(M,k)$ becomes an object in $_H{\cal YD}^H(\alpha^{-1},   
\alpha\beta^{-1}\alpha^{-1})$, with $(h\cdot f) (m)=
f((\alpha^{-1}\beta^{-1}S^{-1}(h))\cdot m)$ and  
$f_{(0)}(m)\ot f_{(1)}=f(m_{(0)})\ot S(m_{(1)})$.  
Moreover, the maps $b_M:k\rightarrow $$\;^*M\ot M$,    
$b_M(1)=\sum_ie^i\ot e_i$ and 
$d_M:M\ot $$\;^*M\rightarrow k$, $d_M(m\ot f)=f(m)$, are morphisms in  
$_H{\cal YD}^H$ and we have  
$(d_M\ot id_M)(id_M\ot b_M)=id_M$ and 
$(id_{\;^*M}\ot d_M )(b_M\ot id_{\;^*M})=id_{\;^*M}$.
\end{proposition}
Consequently, if we consider ${\cal YD}(H)_{fd}$, the subcategory of 
${\cal YD}(H)$ consisting of finite dimensional objects, we obtain:
\begin{theorem} ${\cal YD}(H)_{fd}$ is a braided T-category with left and 
right dualities over $G$, the left (respectively right) duals being given 
as in Proposition \ref{abp} (respectively Proposition \ref{rightdual}).
\end{theorem} 
Assume now that $H$ is finite dimensional. We will construct a 
quasitriangular T-coalgebra over $G$, denoted by $DT(H)$, with the 
property that the T-category $Rep (DT(H))$ of representations of $DT(H)$ 
is isomorphic to ${\cal YD}(H)$ as braided T-categories.\\[2mm]
For $(\alpha , \beta )\in G$, the $(\alpha , \beta )$-component 
$DT(H)_{(\alpha , \beta )}$ will be the diagonal crossed product algebra  
$H^*\bowtie H(\alpha , \beta )$. Define 
\begin{eqnarray*}
&&\Delta _{(\alpha , \beta ), (\gamma , \delta )}:
H^*\bowtie H((\alpha , \beta )*(\gamma , \delta ))\rightarrow 
(H^*\bowtie H(\alpha , \beta ))\ot (H^*\bowtie H(\gamma , \delta )), \\
&&\Delta _{(\alpha , \beta ), (\gamma , \delta )}(p\bowtie h)= 
(p_2\bowtie \gamma (h_1))\ot (p_1\bowtie \gamma ^{-1}\beta \gamma (h_2)).
\end{eqnarray*}
One can check, by direct computation, that these maps are algebra maps, 
satisfying the necessary coassociativity conditions.\\
The counit $\varepsilon $ is just the counit of $DT(H)_{(id, id)}=D(H)$, 
the Drinfeld double of $H$. \\
For $(\alpha , \beta ), (\gamma , \delta )\in G$, define now 
\begin{eqnarray*}
&&\varphi _{(\alpha , \beta )}^{(\gamma , \delta )}:
H^*\bowtie H(\gamma , \delta )\rightarrow 
H^*\bowtie H((\alpha , \beta )*(\gamma , \delta )*(\alpha , \beta )^{-1}), \\
&&\varphi _{(\alpha , \beta )}^{(\gamma , \delta )}(p\bowtie h)=
p\circ \beta \alpha ^{-1}\bowtie \alpha \gamma ^{-1}\beta ^{-1}\gamma (h).
\end{eqnarray*}
Then one can check by direct computation that these are algebra isomorphisms  
giving a {\it conjugation} (that is they are multiplicative and compatible 
with the comultiplications and the counit).\\ 
The antipode is given, for $(\alpha , \beta )\in G$, by 
\begin{eqnarray*}
&&S_{(\alpha , \beta )}:H^*\bowtie H(\alpha , \beta )\rightarrow 
H^*\bowtie H((\alpha , \beta )^{-1}),\\
&&S_{(\alpha , \beta )}(p\bowtie h)=(\varepsilon \bowtie \alpha \beta (S(h))) 
\cdot (S^{*-1}(p)\bowtie 1),
\end{eqnarray*}
where the multiplication $\cdot $ in the right hand side is made 
in $H^*\bowtie H((\alpha , \beta )^{-1})$.\\
Finally, the universal $R$-matrix is given by 
\begin{eqnarray*}
&&R_{(\alpha , \beta ), (\gamma , \delta )}=\sum _i(\varepsilon \bowtie 
\beta ^{-1}(e_i))\ot (e^i\bowtie 1)\in (H^*\bowtie H(\alpha , \beta ))\ot 
(H^*\bowtie H(\gamma , \delta )), 
\end{eqnarray*}
for all $(\alpha , \beta ), (\gamma , \delta )\in G$, where $\{e_i\}$, 
$\{e^i\}$ are dual bases in $H$ and $H^*$.\\
Thus, we have obtained:
\begin{theorem}
$DT(H)$ is a quasitriangular T-coalgebra over $G$, with structure as above.
\end{theorem}
Moreover, the structure of $DT(H)$ was constructed in such a way that, 
via the isomorphisms $_{H^*\bowtie H(\alpha , \beta )}{\cal M}\simeq $$\;
_H{\cal YD}^H(\alpha , \beta )$ from Section \ref{sec2}, we obtain:
\begin{theorem}
$Rep (DT(H))$ and ${\cal YD}(H)$ are isomorphic as braided T-categories 
over $G$.  
\end{theorem}
\begin{remark}{\em 
From ${\cal YD}(H)$ (respectively $DT(H)$) we can obtain, by pull-back 
along the group morphism $Aut_{Hopf}(H)\rightarrow G$, $\alpha  
\mapsto (\alpha , id)$, a braided T-category (respectively a quasitriangular   
T-coalgebra) over $Aut_{Hopf}(H)$. If $\pi $ is a group together with a 
group morphism $\pi \rightarrow Aut_{Hopf}(H)$, by pull-back along it 
we obtain a braided T-category (respectively a quasitriangular T-coalgebra)  
over $\pi $. Quasitriangular T-coalgebras over such $\pi $ have been 
studied by Virelizier in \cite{virelizier}.} 
\end{remark}       
%%%%%%%%%%%%%%%%%%%%%%%%%%%%%%%%%%%%%%%%%%%%%%%%%%%%%%%%%%%%%%%%%%%%%%%%%
\section{An isomorphism of categories $_H{\cal YD}^H(\alpha , \beta )\simeq 
$$\;_H{\cal YD}^H$ \\
in the presence of a pair in involution}\label{sec4}
%%%%%%%%%%%%%%%%%%%%%%%%%%%%%%%%%%%%%%%%%%%%%%%%%%%%%%%%%%%%%%%%%%%%%%%%%%
\setcounter{equation}{0}
%%%%%%%%%%%%%%%%%%%%%%%%%%%%%%%%%%%%%%%%%%%%%%%%%%%%%%%%%%%%%%
The aim of this section is to prove the following result. 
\begin{theorem}
Let $\alpha , \beta $$\;\in Aut_{Hopf}(H)$ and assume that 
there exists $(f, g)$ a  
pair in involution corresponding to $(\alpha , \beta )$. Then the categories  
$_H{\cal YD}^H(\alpha , \beta )$ and $_H{\cal YD}^H$ are isomorphic. \\
A pair of inverse functors $(F, G)$ is given as follows. If 
$M\in $$\;_H{\cal YD}^H(\alpha , \beta )$, then  
$F(M)\in $$\;_H{\cal YD}^H$, where $F(M)=M$ as vector space, with structures 
\begin{eqnarray*}
&&h\rightarrow m=f(\beta ^{-1}(S(h_1)))\beta ^{-1}(h_2)\cdot m,\\
&&m\mapsto m_{<0>}\ot m_{<1>}=m_{(0)}\ot m_{(1)}g^{-1}.
\end{eqnarray*}
If $N\in $$\;_H{\cal YD}^H$, then $G(N)\in $$\;_H{\cal YD}^H
(\alpha , \beta )$,  
where $G(N)=N$ as vector space, with structures
\begin{eqnarray*}
&&h\rightharpoonup n=f(h_1)\beta (h_2)\cdot n, \\
&&n\mapsto n^{(0)}\ot n^{(1)}=n_{(0)}\ot n_{(1)}g. 
\end{eqnarray*}
Both $F$ and $G$ act as identities on morphisms. 
\end{theorem}
\begin{proof}
One checks, by direct computation, that $F$ and $G$ are functors, 
inverse to each other. \\
Alternatively, we can give a very short proof using results from the 
previous section.  
By Example \ref{mp1}, we have  
$_fk^g\in $$\;_H{\cal YD}^H(\alpha , \beta)$.   
By Proposition \ref{abp}, we get    
$(_fk^g)^*\in $$\;_H{\cal YD}^H((\alpha , \beta )^{-1})$. Then, one can check  
that actually $F(M)=(_fk^g)^*\ot M\in \;_H{\cal YD}^H$ and  
$G(N)=(_fk^g)\ot N\in \;_H{\cal YD}^H(\alpha , \beta )$. Also, one can see  
that $(_fk^g)^* \ot $$\;_fk^g=$$\;_fk^g \ot(_fk^g)^*=k$ as objects in 
$_H{\cal YD}^H$, hence $F\circ G=G\circ F=id$,  
using the associativity of the tensor product.
\end{proof}
As we have noticed before, for any $\alpha \;\in Aut_{Hopf}(H)$ we have that 
$(\varepsilon, 1)$ is a pair in involution corresponding to 
$(\alpha , \alpha )$, hence we obtain:
\begin{corollary}
$_H{\cal YD}^H(\alpha , \alpha )\simeq $$\;_H{\cal YD}^H$.
\end{corollary}
Also, as a consequence of the theorem, we obtain the following 
result (a right-left version was given in \cite{doru}), which  
might be useful for the aria of applicability of anti-Yetter-Drinfeld 
modules:    
\begin{corollary}
Assume that there exists a pair in involution $(f, g)$ corresponding to 
$(S^2, id)$. Then the category $_H{\cal YD}^H(S^2, id)$ of   
anti-Yetter-Drinfeld modules is isomorphic to $_H{\cal YD}^H$, and any 
anti-Yetter-Drinfeld module can be written as a tensor product 
$_fk^g\ot N$, with $N\in \;_H{\cal YD}^H$.  
\end{corollary} 
Let again $\alpha ,\beta $$\;\in Aut_{Hopf}(H)$ such that 
there exists $(f, g)$ a   
pair in involution corresponding to $(\alpha , \beta )$, 
and assume that $H$ is 
finite dimensional. Then we know that 
$_H{\cal YD}^H(\alpha , \beta )\simeq    
$$\;_{H^*\bowtie H(\alpha , \beta )}{\cal M}$,   
$_H{\cal YD}^H\simeq $$\;_{D(H)}{\cal M}$, and the isomorphism 
$_H{\cal YD}^H(\alpha , \beta )\simeq $$\;_H{\cal YD}^H$ constructed in the  
theorem is induced by an algebra isomorphism between 
$H^*\bowtie H(\alpha , \beta )$ and $D(H)$, given by  
\begin{eqnarray*}
&&D(H)\rightarrow H^*\bowtie H(\alpha , \beta ),\;\;\;
p\ot h\mapsto g^{-1}\rightharpoonup p\bowtie f(\beta ^{-1}(S(h_1)))
\beta ^{-1}(h_2), \\
&&H^*\bowtie H(\alpha , \beta )\rightarrow D(H), \;\;\;
p\bowtie h\mapsto g\rightharpoonup p\ot f(h_1)\beta (h_2).
\end{eqnarray*}       
%%%%%%%%%%%%%%%%%%%%%%%%%%%%%%%%%%%%%%%%%%%%%%%%%%%%%%%%%%%%


\begin{thebibliography}{99}
%%%%%%%%%%%%%%%%%%%%%%%%%%%%%%%%%%%%%%%%%%%%%%%%%%%%%%%%%%%%

\bibitem{bpvo}
D. Bulacu, F. Panaite, F. Van Oystaeyen, Generalized diagonal
crossed products and smash products for (quasi) Hopf algebras. Applications,
in preparation.

\bibitem{stef}
S. Caenepeel, M. De Lombaerde, A categorical approach to Turaev's 
Hopf group-coalgebras, arXiv:math.QA/0409600. 

\bibitem{cmz}
S. Caenepeel, G. Militaru, S. Zhu, ''Frobenius and separable
functors for generalized module categories and nonlinear
equations'', {\sl Lecture Notes in Mathematics} {\bf 1787},
Springer Verlag, Berlin, 2002.

\bibitem{cvoz}
S. Caenepeel, F. Van Oystaeyen, Y. Zhang, The Brauer group of Yetter-Drinfeld 
module algebras, {\sl Trans. Amer. Math. Soc.} {\bf 349} (1997), 
3737--3771.

\bibitem{connes1}
A. Connes, H. Moscovici, Hopf algebras, cyclic cohomology and the 
transverse index theorem, {\sl Comm. Math. Phys.} {\bf 198} (1998), 
199--264. 

\bibitem{connes2}
A. Connes, H. Moscovici, Cyclic cohomology and Hopf algebra symmetry, 
{\sl Lett. Math. Phys.} {\bf 52} (2000), 1--28.  

\bibitem{hajac1}
P. M. Hajac, M. Khalkhali, B. Rangipour, Y. Sommerhauser,  
Stable anti-Yetter-Drinfeld modules, {\sl C. R. Math. Acad. Sci. Paris}  
{\bf 338} (2004), 587--590.

\bibitem{hajac2}
P. M. Hajac, M. Khalkhali, B. Rangipour, Y. Sommerhauser, 
Hopf-cyclic homology and cohomology with coefficients, {\sl C. R. Math. 
Acad. Sci. Paris} {\bf 338} (2004), 667--672.     

\bibitem{hn1}
F. Hausser and F. Nill, Diagonal crossed products by duals of
quasi-quantum groups, {\sl Rev. Math. Phys.} {\bf 11} (1999),
553--629.

\bibitem{jara}
P. Jara, D. \c Stefan, Cyclic homology of Hopf Galois extensions and 
Hopf algebras, arXiv:math.KT/0307099. 

\bibitem{k}
C. Kassel, ''Quantum groups'', {\sl Graduate Texts in Mathematics}
{\bf 155}, Springer Verlag, Berlin, 1995.

\bibitem{m}
S. Majid, ''Foundations of quantum group theory'', Cambridge Univ.
Press, 1995.

\bibitem{mon}
S. Montgomery, ''Hopf algebras and their actions on rings'', CBMS Regional 
Conference Series, Vol. 82, Amer. Math. Soc., Providence, RI, 1993.

\bibitem{doru}
M. D. Staic, A note on anti-Yetter-Drinfeld modules, preprint 2004.

\bibitem{sw}
M. E. Sweedler, ''Hopf algebras'', Benjamin, New York, 1969.

\bibitem{turaev}
V. Turaev, Homotopy field theory in dimension 3 and crossed group-categories, 
arXiv:math.GT/0005291. 

\bibitem{virelizier}
A. Virelizier, Graded quantum groups, arXiv:math.QA/0312330.

\bibitem{zunino}
M. Zunino, Yetter-Drinfeld modules for crossed structures, 
{\sl J. Pure Appl. Algebra} {\bf 193} (2004), 313--343. 

\end{thebibliography}
\end{document}